\documentclass[11pt]{article}

\usepackage{amscd,amsmath, amssymb, fancyhdr, url,epsfig, color, enumitem, tocloft}
\usepackage{dutchcal}
\usepackage{graphicx}

\usepackage[backref=page]{hyperref}
\renewcommand*{\backrefalt}[4]{%
	\ifcase #1 (Not cited.)%
	\or        (Cited on page~#2.)%
	\else      (Cited on pages~#2.)%
	\fi}

\hypersetup{
	colorlinks   = true,
	citecolor    = magenta,
	linkcolor    =blue,
	urlcolor     =magenta	
}

\newcommand{\version}{version 2.2,\ \ March 16, 2021}
\makeatletter
\def\x@arrow{\DOTSB\Relbar}
\def\xlongrightarrowfill@{\arrowfill@\relbar\relbar\longrightarrow}
\newcommand{\xlongrightarrow}[2][]{%
        \ext@arrow 0099\xlongrightarrowfill@{#1}{#2}}
\makeatother

\newcommand{\la}{\lambda}
\newcommand{\al}{\alpha}

\newcommand{\f}{\varphi}

\newcommand{\Lie}{\operatorname{Lie}}

\numberwithin{equation}{section}

\def\eqref#1{(\ref{#1})}
\newcommand{\goth}{\mathfrak}

\newcommand{\ra}{{\:\longrightarrow\:}}
\newcommand{\Z}{{\mathbb Z}}
\newcommand{\C}{{\mathbb C}}
\newcommand{\R}{{\mathbb R}}
\newcommand{\Q}{{\mathbb Q}}

\newcommand{\6}{\partial}
\def\1{\sqrt{-1}\:}
\newcommand{\restrict}[1]{{\left|_{{\phantom{|}\!\!}_{#1}}\right.}}
\newcommand{\cntrct}                
{\hspace{2pt}\raisebox{1pt}{\text{$\lrcorner$}}\hspace{2pt}}
\newcommand{\arrow}{{\:\longrightarrow\:}}

\renewcommand{\bar}{\overline}
\renewcommand{\phi}{\varphi}
\renewcommand{\epsilon}{\varepsilon}
\renewcommand{\geq}{\geqslant}

\newcommand{\End}{\operatorname{End}}
\newcommand{\Tot}{\operatorname{Tot}}

\newcommand{\Aut}{\operatorname{Aut}}

\newcommand{\Diff}{\operatorname{Diff}}

\newcommand{\codim}{\operatorname{codim}}

\newcommand{\GL}{\operatorname{GL}}

\newcounter{Mycounter}[section]
\newcounter{lemma}[section]
\newcounter{claim}[section]
\renewcommand{\theclaim}{{Claim \thesection.\arabic{claim}}}
\newcommand{\claim}{%
     \setcounter{claim}{\value{Mycounter}}
     \refstepcounter{claim}
     \stepcounter{Mycounter}
     {\noindent \bf \theclaim:\ }}

\newcounter{sublemma}[section]
\newcounter{corollary}[section]
\renewcommand{\thecorollary}{{Corollary \thesection.\arabic{corollary}}}
\newcommand{\corollary}{%
     \setcounter{corollary}{\value{Mycounter}}
     \refstepcounter{corollary}
     \stepcounter{Mycounter}
     {\noindent \bf \thecorollary:\ }}

\newcounter{theorem}[section]
\renewcommand{\thetheorem}{{Theorem \thesection.\arabic{theorem}}}
\newcommand{\theorem}{%
     \setcounter{theorem}{\value{Mycounter}}
     \refstepcounter{theorem}
     \stepcounter{Mycounter}
     {\noindent \bf \thetheorem:\ }}

\newcounter{conjecture}[section]
\newcounter{proposition}[section]
\renewcommand{\theproposition} {{Proposition \thesection.\arabic{proposition}}}
\newcommand{\proposition}{%
     \setcounter{proposition}{\value{Mycounter}}
     \refstepcounter{proposition}
     \stepcounter{Mycounter}
     {\noindent \bf \theproposition:\ }}

\newcounter{definition}[section]
\renewcommand{\thedefinition} {{Definition~\thesection.\arabic{definition}}}
\newcommand{\definition}{%
     \setcounter{definition}{\value{Mycounter}}
     \refstepcounter{definition}
     \stepcounter{Mycounter}
     {\noindent \bf \thedefinition:\ }}

\newcounter{example}[section]
\renewcommand{\theexample}{{Example \thesection.\arabic{example}}}
\newcommand{\example}{%
     \setcounter{example}{\value{Mycounter}}
     \refstepcounter{example}
     \stepcounter{Mycounter}
     {\noindent \bf \theexample:\ }}

\newcounter{remark}[section]
\renewcommand{\theremark}{{Remark \thesection.\arabic{remark}}}
\newcommand{\remark}{%
     \setcounter{remark}{\value{Mycounter}}
     \refstepcounter{remark}
     \stepcounter{Mycounter}
     {\noindent \bf \theremark:\ }}

\newcounter{problem}[section]
\newcounter{question}[section]
\makeatletter

\@addtoreset{equation}{section}
\@addtoreset{footnote}{section}
\makeatother

\def\blacksquare{\hbox{\vrule width 5pt height 5pt depth 0pt}}
\def\endproof{\blacksquare}

\newcommand{\proof}{{\bf Proof: \ }}

\begin{document}

\begin{center}
{\Large\bf  Closed orbits of Reeb fields on Sasakian manifolds and elliptic curves on Vaisman manifolds}\\[5mm]

{\large
Liviu Ornea\footnote{Liviu Ornea is  partially supported by UEFISCDI grant PCE 30\slash 2021. },  
Misha
Verbitsky\footnote{Misha Verbitsky is partially supported by the 
Russian Academic Excellence Project '5-100'', FAPERJ E-26/202.912/2018 
and CNPq - Process 313608/2017-2.\\[1mm]

}\\[4mm]

}

\end{center}

{\small
\hspace{0.15\linewidth}
\begin{minipage}[t]{0.72\linewidth}
{\bf Abstract} \\ 
A compact complex manifold $V$ is called Vaisman if
it admits an Hermitian metric which is
conformal to a K\"ahler one, and 
a non-isometric conformal action by
$\Bbb C$. It is called quasi-regular if the  
$\Bbb C$-action has closed orbits.
In this case the corresponding
leaf space is a projective orbifold,
called the quasi-regular quotient of $V$.
It is known that the set of all quasi-regular
Vaisman complex structures is dense in the 
appropriate deformation space.
We count the number of closed elliptic
curves on a Vaisman manifold, proving
that their number is either infinite
or equal to the sum of all Betti
numbers of a K\"ahler orbifold obtained as
a quasi-regular quotient of $V$.
We also give a new proof of a result by
Rukimbira showing that the number of 
Reeb orbits on a Sasakian manifold $M$ 
is either infinite or equal to the sum of all Betti
numbers of a K\"ahler orbifold obtained as
an $S^1$-quotient of $M$.
\end{minipage}
}

\tableofcontents


\section{Introduction}
\label{_Intro_Section_}

Initially stated for periodic orbits of Hamiltonian flows
on hypersurfaces of contact type in symplectic manifolds,
Weinstein conjecture (\cite{_Weinstein_}) can be
intrinsically formulated in the following form:

\hfill

\noindent{\bf Weinstein Conjecture:} On any closed contact
manifold $(N,\eta)$ the Reeb field has at least one closed
orbit.

\hfill

For history and background on Weinstein conjecture, 
see the survey \cite{_Pasquotto:W_conj_}.

\hfill

The conjecture has been solved in the affirmative  for
contact hypersurfaces in Euclidean space by C. Viterbo
(\cite{_Viterbo_}), on $S^3$, with a different method, by
H. Hofer (\cite{_Hofer_}) and, more generally, on closed
3-manifolds, by C.H. Taubes (\cite{_Taubes_}). Several
extensions of these results are also available. But in
full generality, the conjecture is still open and remains
one of the most challenging problems in symplectic and
contact topology.

The dynamic of the Reeb flow on a contact manifold
became a separate subfield in symplectic geometry, called 
the Reeb dynamics (\cite{_Geiges:lectures_}). One of the central questions of this
subfield is finding an estimate for the number of closed 
Reeb orbits on contact manifolds under different
geometric conditions.

Most research in this direction is based on 
contact homology, which is a Morse-type cohomology theory
expressing the Reeb orbits as fixed point of a certain
gradient acrion on an appropriate loop space. As in
Morse (and Floer) theory, the number of closed Reeb orbits bounds
the sum of Betti numbers of the relevant cohomology
theory.

The most striking achievement in this direction 
is a result of Hutchins and Cristofaro-Gardiner,
who proved that any 3-dimensional compact contact
manifolds admits at least 2 distinct closed Reeb
orbits (\cite{_Cristofaro-Gardiner_Hutchins_}). 

However, for dimension $> 3$ this approach fails,
and in general nothing is known.

Adding geometric structures compatible with the given
contact structure puts an interesting twist on the
Reeb dynamics. In this paper we study the Reeb
dynamics on Sasakian manifolds.

 These are contact manifolds equipped with a
Riemannian metric, in such a way that
their symplectic cone is equipped with a K\"ahler structure which is
automorphic with respect to the dilations along the
generators of the cone (see Subsection \ref{_contact_subsection_}). In
this sense, Sasakian manifolds can be understood as
odd-dimensional counterparts of K\"ahler manifolds. They
are ubiquitous and their geometry is a well established
subject, see \cite{_Boyer_Galicki_}.

The first main result of this paper is a new proof of a result in \cite{_Rukimbira_} which confirms 
Weinstein conjecture for compact Sasakian manifolds:

\hfill

\theorem\label{_Main_} Let $M$ be a compact Sasakian
manifold of dimension $2n+1$. Then its Reeb field has at
least $n+1$ closed orbits.

\hfill

For the proof, we first observe that the Reeb field can be
approximated by a quasi-regular one, i.e. with a Reeb
field  such that its space of orbits is a projective
orbifold. We then relate the closed orbits of the initial
Reeb field to the fixed points of the action of the group
generated by its flow on this projective
orbifold. Finally, we apply the orbifold version,
\cite{_Fontanari_}, of a celebrated result of A. Bia\l
ynicki-Birula, \cite{_Birula_}, counting fixed points of
algebraic groups.

Our proof does not explicitly use the contact
or Sasakian geometry, but transfers the problem to the
framework of complex geometry.

From this argument we could also obtain an explicit
expression for the number of the Reeb orbits on a Sasakian
manifold associated
with a projective orbifold $X$ admitting  holomorphic
vector fields with isolated zeros. In that case
the number of  Reeb orbits (for an appropriate Sasakian
structure) is equal to the number of fixed points
of any of these vector fields. It is also equal to the
sum of Betti numbers of $S$. This always gives the
lower estimate, and gives an upper estimate
when the Sasakian structure is sufficiently
general.

The compact Sasakian manifolds are known to be closely
related  to Vaisman manifolds (see Section \ref{_vaisman_section_}
for their definition and properties). Indeed, Vaisman manifolds are mapping
tori of circles with fibres compact Sasakian manifolds. Diagonal Hopf
manifolds and their complex compact submanifolds are
typical examples. 

We use a similar approach to count the 
elliptic curves on compact Vaisman manifolds. The result is not a corollary of the one for Sasakian manifolds, but the idea of the proof is similar. Our second main result reads:

\hfill

\theorem\label{_elliptic_curves_on_vaisman_} 
Let $V$ be a compact Vaisman manifold which contains precisely $r$ elliptic curves. Then
$r\geq \dim_\C V$.

\proof \ref{_Vaisman_bound_from_Lfsch_Corollary_}. \endproof

\hfill

A more precise result 
is stated in \ref{_Vaisman_number_Theorem_}.

\hfill

This theorem is proven in Section \ref{_vaisman_section_}.
Note that this result is elementary for the diagonal Hopf
manifolds. 

\hfill

After the first version of  this paper was finished,
we became aware of the papers (\cite{_Rukimbira_}, \cite{_Goertsches_Nozawa_Tobin_}), 
in particular the work of  Philippe Rukimbira
who obtained the same bounds on the number of 
closed Reeb orbits using different methods, suitable in the non-integrable case. Our method is specific to complex geometry and 
the result about the number of 
elliptic curves on Vaisman manifolds
is new. We are grateful to Gianluca Bande 
who first sent us the reference and to Oliver Goertsches
and Hiraku Nozawa for clarifying to us a part of their work.


\section{CR, contact and Sasakian geometries}


We present the necessary notions concerning Sasakian
geometry. For details and examples, see
\cite{_Boyer_Galicki_}. We start by recalling the
definition of CR and contact structures. Both of these
geometric structures are subjancent to the Sasakian structures.

\subsection{CR manifolds}

\definition 
	Let $M$ be a smooth manifold,
	$B\subset TM$ a sub-bundle in the tangent bundle,
	and $I:\; B \arrow B$ an endomorphism satisfying
	$I^2=-1$. Consider its $\1$-eigenspace $B^{1,0}(M)\subset
	B\otimes \C \subset T_C M=TM\otimes \C$.
	Suppose that $[B^{1,0}, B^{1,0}]\subset B^{1,0}$.
	Then $(B,I)$ is called a {\bf CR-structure on $M$}.

\hfill

\definition Let $(M,B)$ be a CR manifold and  $\Pi_{TM/B}:TM\ra TM/B$ be 
	the projection 
	to the normal bundle of $B$ in $TM$. The  tensor field $B\otimes B \arrow TM/B$ mapping
	vector fields $X, Y\in B$ 
	to $\Pi_{TM/B}([X,Y])$ is called {\bf the Frobenius form
		of $B$}. It is the obstruction 
	to the integrability of the distribution given by $B$.

\hfill

\remark Let $S$ be a CR manifold, with  the  bundle $B$
of codimension 1, and almost complex structure
$I\in \End(B)$. Since the Frobenius form vanishes 
when both arguments are from  $B^{0,1}$ and $B^{1,0}$, 
it is a pairing between $B^{0,1}$ and $B^{1,0}$. Indeed, $[B^{1,0}, B^{1,0}]\subset B^{1,0}$
and $[B^{0,1}, B^{0,1}]\subset B^{0,1}$. This proves that
the Frobenius form is a Hermitian form taking values in a
trivial rank 1 bundle $TM/B$. The Frobenius form on a CR
manifold $(M,B,I)$
with $\codim B=1$ is called {\bf the Levi form}.

\hfill

\remark
If, in addition, $B\subset TM$ is a contact bundle,
its Levi form is non-degenerate. Therefore,
it has constant rank. If it is positive or negative
definite, the CR manifold $(M,B,I)$ is called
{\bf strictly pseudoconvex}. In this case
we fix the orientation on the trivial bundle 
$TM/B$ in such a way that the Levi form is
positive definite.

\hfill

\example \begin{enumerate}
	\item[(i)] A complex manifold $(X,I)$ is CR, with $B=TX$. Indeed $[T^{1,0}X,T^{1,0}X]\subset T^{1,0}X$ is equivalent to the Newlander-Nirenberg theorem.
	\item[(ii)] Let $(X,I)$ be a complex manifold and $M\subset X$ a real hypersurface. Then 
	$B:=TX\cap I(TX)$ is a distribution of dimension $\dim_\C X-1$ which gives $M$ the structure of a CR manifold.
	\item[(iii)] Let $X$ be a complex manifold and $\f:X\to\R$ a strictly plurisubharmonic function. Then $\omega:=\6\bar\6\f$ is positive definite and all level sets $M_c:=\f^{-1}(c)$ are strictly pseudoconvex CR manifolds, with Frobenius forms $\omega\restrict{M_c}$.
	\end{enumerate}

The following result was proved by D. Burns, but never published by him, see \cite{_Lee_}; another proof, valid also for $\dim M=3$, is given in \cite{_Schoen_}.

\hfill

\theorem {\bf (D. Burns)}\label{_Compact_CR_group_}
Suppose $M$ is a compact, connected, strictly pseudoconvex CR
manifold of dimension $2n+1\geq 5$. The full CR automorphism group $\Aut(M,B,I)$ is compact unless $M$
is globally CR equivalent to $S^{2n+1}$ with its standard CR structure.

\hfill

\remark CR automorphisms are also called {\bf CR
  holomorphic diffeomorphisms}. A vector field whose flow 
consists in CR automorphisms is called {\bf a CR holomorphic vector field.}

\subsection{Contact manifolds}\label{_contact_subsection_}

\definition Let $S$ be a manifold. Then $C(S):= S
\times \R^{>0}$ is called {\bf the  cone over $S$}. The
multiplicative group $\R^{>0}$ acts on $C(S)$ by dilations along the generators: 
$h_\la(x, t) \arrow (x, \lambda t)$.

\hfill

\definition A differential $k$-form $\alpha$ on a  cone is called {\bf  automorphic} if $h_q^*\al=q^k\al$.

\hfill

\definition A {\bf contact manifold} is 
a manifold $S$ such that its cone $C(S)$ is endowed with an automorphic 
symplectic form $\omega$, i.e. $h_\la^*\omega=\la^2\omega$. 

\hfill

The next characterization explains the geometry of a contact manifold and its relation to CR geometry:

\hfill

\theorem Let $S$ be a differentiable manifold. The following are equivalent:
\begin{enumerate}
\item[(i)]  $S$ is contact.
\item[(ii)] $S$ is odd-dimensional and there exists an oriented sub-bundle of codimension 1,
$B\subset TS$,  
with non-degenerate Frobenius form $\Lambda^2 B \stackrel \Phi \arrow TS/B$. 
This $B$ is called  {\bf the contact bundle}.

\item[(iii)] $S$ is odd-dimensional and there exists an oriented sub-bundle of codimension 1,
$B\subset TS$, such that for any nowhere degenerate 1-form  $\eta\in
\Lambda^1 S$ which annihilates $B$, the form $\eta \wedge (d\eta)^k$ 
is a  volume form (where $\dim S=2k+1$). Every such $\eta$ 
is called {\bf a contact form}.
\end{enumerate}

\proof Well known (see \cite{_McDuff_Salamon_}).
\endproof

\hfill

\definition 	Let $S$ be a contact manifold with a  contact form $\eta$. The vector field $R$ defined by the equations:
$$R\cntrct \eta=1,\qquad R\cntrct d\eta=0.$$
is called {\bf the Reeb} or {\bf characteristic} field.

\subsection{Sasakian manifolds}

\definition The {\bf Riemannian cone} of a Riemannian manifold $(M,g_M)$ is 
$C(M):=M\times \R^{>0}$ endowed with the metric $g=t^2g_M+dt\otimes dt$, where $t$ is the coordinate on $\R^{>0}$.

\hfill

\definition A {\bf Sasakian structure} on a Riemannian
manifold $(M,g_M)$ is a K\"ahler structure $(g,I,\omega)$ on its
Riemannian cone  $(C(M), g=t^2g_M+dt\otimes dt)$  such that:
\begin{itemize}
	\item $I$ is $h_\la$ - invariant, and
	\item $\omega$ is automorphic: $h_\la^*\omega=\la^2\omega$.
\end{itemize}
In this case, $(C(M), g, I, \omega)$ is called {\bf the cone} of the Sasakian manifold.

\hfill

\remark \label{_Sas_properties_}
\begin{enumerate}
	\item[(i)] It can be seen that $\f=t^2$ is a
          K\"ahler potential for $\omega$, that is,
          $\omega = dd^c \f$, where 
$d^c = I d I^{-1}$. The manifold  $M$ can be
          identified with a level of $\f$. 
The contact bundle on $M$ can be obtained
as $TM \cap I(TM)$, where $TM$ is considered as
a sub-bundle of $TC(M)$, and $I$ the complex structure on
$TC(M)$. This gives a CR-structure on $M$.
An easy calculation implies that the restriction of 
$dd^c \f$ to the contact bundle is equal to the Levi form.
Therefore, the CR-structure on the Sasakian manifolds
is strictly pseudoconvex.

	\item[(ii)] Let $\xi:=t\frac{t}{dt}$ be the {\bf Euler field} on $C(M)$. The complex structure of the cone is $h_\la$ invariant, thus the Euler field acts on the cone by holomorphic homotheties.
	\item[(iii)] Since $\omega$ is, in particular, a
          symplectic form, a Sasakian manifold is, in
          particular, a {\bf contact manifold}. Its
          contact form is defined as follows. Let
          $\eta:=\xi\cntrct\omega$. Then, by Cartan's
          formula, $d\eta=\Lie_\xi\omega=\omega$, and
          hence
          $(d\eta)^n\wedge\eta=\frac{1}{n+1}\xi\cntrct\omega^{n+1}$
          is non-degenerate, thus each slice
          $M\times\{t_0\}\subset C(M)$ is contact, with
          contact form the restriction of $\eta$. We may
          consider $\eta$ as a 1-form on $M$. 
	\item[(iv)] The above CR structure is given by the distribution $B:=\ker\eta$.
	\item[(v)] Let $R:=I\xi$. Clearly, $R$ is tangent to $M$ and transverse to the CR distribution. Then $\eta(R)=1$ and $R\cntrct d\eta=0$, and hence $R$ is the {\bf Reeb field} of the contact manifold $(M,\eta)$. Note that on a Sasakian manifold, the Reeb field is the $g_M$-dual of the contact form: $R=\eta^\sharp$. Moreover, its flow consists of contact isometries: $\Lie_R g_M=0$, $\Lie_R\eta=0$.
\end{enumerate}

The relation between CR structures and Sasakian structures is clarified in the following result:

\hfill

\theorem {\bf (\cite{_Ornea_Verbitsky_CR_})}\label{_CR_versus_Sas_}  Let  $(B,I)$ be a CR structure on an odd-dimensional compact manifold $M$. Then there exists a compatible Sasakian  metric on $M$ if and only if $M$ admits a CR-holomorphic vector field which is transversal to $B$. Moreover, for every such field $v$, there exists a unique Sasakian metric such that $v$ or $-v$ is its Reeb field.

\hfill

\example Odd-dimensional spheres with the round metric
 are equipped with a natural Sasakian structure, since their
cones are $\C^n\setminus 0$ on which the standard flat K\"ahler
form $\sqrt{-1}dz_i\wedge d\bar z_i$ is automorphic.

\hfill

\definition A Sasakian manifold is called {\bf quasi-regular} if all orbits of its Reeb field are compact.

\hfill

If $M$ is a compact quasi-regular sSsakian manifold, with
Reeb field $R$, one can consider the space of orbits
$X:=M/R$ which is the same as $C(M)/\langle \xi,
R\rangle$. In this case $C(M)$ is the total space
of a principal holomorphic $\C^*$-bundle over $X$ and the
corresponding line bundle has positive curvature
$dd^c\log\f$. Therefore, Kodaira's theorem implies
(\cite{_Ornea_Verbitsky_immersion_}):

\hfill

\theorem 
Let $M$ be a compact, quasi-regular Sasakian
manifold, with Reeb field $R$. Then the space of orbits
$X:=M/R$ is a projective orbifold.

\hfill

\example \begin{enumerate}
	\item Let $X\subset \C P^n$ be a complex submanifold, and $C(X)\subset  \C^{n+1}\setminus 0$ the
corresponding cone. The cone $C(X)$ is obviously K\"ahler and its K\"ahler form is automorphic,
hence the intersection $C(X)\cap S^{2n-1}$ is Sasakian. This intersection is an $S^1$-
bundle over $X$. This construction gives many interesting contact manifolds,
including Milnor's exotic 7-spheres, which happen to be Sasakian.
\item In general, every link of homogeneous singularity is Sasakian. All quasi-regular Sasakian manifolds are obtained this way.

\item All 3-dimensional Sasakian manifolds are quasi-regular, see \cite{_Belgun_, _Geiges_}.
\end{enumerate}

The following result shows that on a compact manifold, a Sasakian structure can be aproximated with quasi-regular ones.

\hfill

\theorem {\bf
  (\cite{_Ornea_Verbitsky_immersion_})}\label{_qr_approximation_th_}
Let $M$ be a compact Sasakian manifold with Reeb field
$R$. Then $R=\lim_i R_i$, where $R_i$ are the Reeb fields
of quasi-regular Sasakian structures on $M$.

\hfill

\proof We present a proof slightly different from the
original one, without making use of \ref{_Compact_CR_group_}. Let $(B, I)$ be the subjacent CR structure,
which is strictly pseudoconvex, \ref{_Sas_properties_},
(i). Let $G$ be the closure of the Lie subgroup generated
by the flow of $R$ in the group of CR automorphisms
$\mathrm{Aut}(M,B,I)$. Since $R$ acts on $M$ by
isometries,  $G$ is compact. Being also
commutative, $G$ is a torus. Now, any vector field $R'$
sitting in the Lie algebra of $G$ sufficiently close to
$R$ will still be transversal to $B$, and hence it will be
the Reeb field of another Sasakian structure
(\ref{_CR_versus_Sas_}).

But a Reeb field is quasi-regular if and only if it generates a compact subgroup, i.e. it is rational with respect to the rational structure of the Lie algebra of $G$. However, the set of rational points is dense in this Lie algebra. \endproof

\hfill

\remark\label{_complex_structure_on_quotient_} It is important to note that during the above approximation process, the CR structure remains unchanged. In particular, for each quasi-regular Reeb field $R_i$ approximating $R$, the complex structure on the projective manifold $M/R_i$ is the same, and corresponds to the one on the contact distribution on $M$.


\section{Closed orbits of Reeb fields}


We can prove now the first of our results, \ref{_Main_}.

Let $M$ be a $(2n+1)$-dimensional Sasakian manifold, with Reeb field $R$ and subjacent CR structure $(B,I)$. As in the proof of \ref{_qr_approximation_th_}, denote with $G$ the closure of the group generated by the flow of $R$ in the group of CR automorphisms $\mathrm{Aut}(M,B,I)$, i.e. $G:=\overline{\langle e^{tR}\rangle}$, for $t\in\R$. The following statement is then clear:

\hfill

\claim There exists a one-to-one correspondence between 1-dimensional orbits of $G$ and  closed orbits of the Reeb field $R$. 

\hfill

Let $R'$ be a quasi-regular approximation of $R$ (see \ref{_qr_approximation_th_}) and $X=M/R'$ the corresponding projective quotient  orbifold. The key observation is:

\hfill

\remark The group $G$ acts on $X$ by holomorphic isometries. Moreover, there is a one-to-one correspondence between 1-dimensional orbits of the action of $G$ on $M$ and fixed points of the action of $G$ on $X$.

\hfill

We thus reduced the problem of counting closed orbits of
the Reeb field on $M$ to counting fixed points of a group
acting by holomorphic isometries on a projective
orbifold. This, in turn, is equal to the number of zeros
of a generic vector field $\zeta\in\Lie(G)$. For compact
projective manifolds, this number is computed by a
celebrated theorem of A. Bia\l ynicki-Birula,
\cite{_Birula_}. The orbifold version that we need is due
to E. Fontanari:

\hfill

\theorem \label{_Fontanari_Theorem_}
{(\cite{_Fontanari_})} Let $\zeta$ be a holomorphic
vector field on a compact projective orbifold $X$ of
complex dimension $k$. Then the number of zeros of $\zeta$ is
equal to $\sum_{i=0}^{2k}b_i(X)$, the sum of all Betti numbers of $X$.

\hfill

Observe that by Lefschetz theorem, the sum of all Betti
numbers is at least $\dim_\C X+1$, and in our case
$\dim_\C X=n$, which completes the proof of
\ref{_Main_}. \endproof

\hfill

\remark Note that the number of zeros exhibited in 
Bia\l ynicki-Birula's and Fontanari's theorems do not count
multiplicities, that is, there exist at least $n+1$ {\em
  distinct} closed orbits.


\section{Vaisman manifolds}\label{_vaisman_section_}


\subsection{Definition and the canonical foliation}
Vaisman manifolds are a significant and much studied
subclass of locally conformally K\"ahler manifolds, see
\cite{_Dragomir_}. Here we give a definition which is
suitable for our purpose:

\hfill

\definition \label{_Vaisman:def_}
Let $(V,I,g_V)$ be a Hermitian manifold such
that the fundamental form
$\omega_V(\cdot,\cdot)=g_V(\cdot,I\cdot)$ satisfies the
equation $d\omega_V=\theta\wedge\omega_V$ for a
non-zero, $\nabla^{g_V}$-parallel 1-form $\theta$, where
$\nabla^{g_V}$ is the Levi-Civita connection of
$g_V$. Then $(V,I,g_V)$ is a {\bf Vaisman manifold} and
$\theta$ is its {\bf Lee form}.

\hfill

\remark By \cite[Subsection 1.3]{_Ornea_Verbitsky_LCK_rank_}, all compact
Vaisman manifolds, considered as complex manifolds, are
obtained in the  following way. Let
$C(M)$ be the K\"ahler cone of a compact Sasakian
manifold and $q$ a non-trivial holomorphic homothety of
$C(M)$. Then the compact complex manifold $C(M)/\langle
q\rangle$ is Vaisman. 

\hfill

\example\label{_examples_} Diagonal Hopf manifolds
$(\C^n\setminus 0)/\langle A\rangle$ where $A\in\GL(n,\C)$
is diagonalizable, with eigenvalues of absolute value
strictly greater than $1$, are Vaisman. All compact complex 
submanifolds of a Vaisman manifold are
Vaisman. Non-K\"ahler elliptic surfaces are Vaisman
(\cite{_Belgun_}).

\hfill
 
\remark \begin{enumerate}
	\item[(i)] The 2-form $\omega_0:=dd^c\log t$ on
          $C(M)=M\times \R^{>0}$ is
          $q$-invariant. Moreover, one can see that
          $\omega_0:=\frac{1}{t^2}\left(\omega_V-dt\wedge
          I(dt)\right)$, hence $\omega_0$ is
          positive-definite in the directions transversal to
          $\left\langle \frac d{dt}, I\left(\frac
          {d}{dt}\right)\right\rangle$, 
and  vanishes on $\left\langle \frac d{dt}, I\left(\frac
          {d}{dt}\right)\right\rangle$.
	\item[(ii)] The 1-form $d\log t$, which lives on
          the cone, is already $q$-invariant. Therefore,
          the 2-form $\omega_0$ descends to an exact form
          on the Vaisman manifold $V$.
	\end{enumerate}

\definition \label{_cano_foli_Definition_}
The foliation $\Sigma:=\ker\omega_0$ on $V$ is called {\bf the canonical foliation} of the Vaisman manifold. 

\hfill

\remark \label{_Lee_field_Remark_}
The foliation $\Sigma$ is  generated by the $g_V$-duals of $\theta$ and $I\theta$, which are commuting, Killing and real holomorphic vector fields.

\hfill

The name of this foliation is motivated by the following:

\hfill

\proposition The foliation $\Sigma$ is independent on the
choice of the cone $C(M)$ and of the choice of the
holomorphic homothety $q$.

\hfill

\proof Suppose we have two different exact and (semi-)positive forms $\omega_0$ and $\omega_0'$. Then the sum $\omega_1:=\omega_0+\omega_0'$ is still exact and (semi-)positive. If $\ker\omega_0\neq \ker\omega_0'$, the form $\omega_1$ is strictly positive, which is impossible because $\omega_0$ and $\omega_0'$ are exact and then Stokes theorem implies $\int_V\omega_1^{\dim_\C V}=0$. \endproof

\subsection{Complex curves on Vaisman manifolds are elliptic}

Complex curves are very particular on compact Vaisman
manifolds: they have to be elliptic, as shown by the next
result.

\hfill

\theorem Let $C$ be a complex curve on a compact Vaisman manifold $V$. Then $C$
is a leaf of the canonical foliation. In particular, $C$ is an elliptic curve.

\hfill

\proof As above, $\int_C\omega_0=0$, by Stokes, hence $C$
is tangent to $\Sigma=\ker\omega_0$. But all compact
leaves of $\Sigma$ are elliptic since the tangent bundle
$T\Sigma$ is trivial by construction. \endproof

\hfill

\remark We proved in \cite{_OV:Potential_,_Ornea_Verbitsky_immersion_}
that a compact Vaisman manifold $V$ can be holomorphically
embedded in a diagonal Hopf manifold $H$ (see
\ref{_examples_}). Intersecting $V$ with two complementary
flags of Hopf submanifolds (which exist due to a result of Ma. Kato, 
\cite{_Kato_}), we see that $V$ contains at least two
elliptic curves. In fact, there are many more elliptic
curves, as shown below.


\subsection{Elliptic curves on a Vaisman manifold and
  group orbits}


Let $V$ be a Vaisman manifold. The Vaisman structure,
as given in \ref{_Vaisman:def_}, is a notion of
differential-geometric nature. However,
there exists a complex-geometric counterpart,
which can be defined as follows.
Recall that {\bf the LCK rank} (\cite{_Ornea_Verbitsky_LCK_rank_}) of a Vaisman
manifold is the rank of the smallest rational
subspace in $H^1(V,\Q)$ containing the cohomology
class of the Lee form. The monodromy group of the
smallest K\"ahler covering of $V$ is $\Z^d$,
where $d$ is its LCK rank.
From \cite{_Ornea_Verbitsky_LCK_rank_} it follows that
any complex manifold of Vaisman type admits
a Vaisman metric of LCK rank 1, and, moreover,
such metrics are dense in the set of all
Vaisman metrics. If we are interested
in complex geometry, we may always assume
that $V$ has LCK rank 1, and admits
a K\"ahler $\Z$-cover $\tilde V$.

It is well-known (see, for example,
\cite[Section 5]{_Ornea_Verbitsky_immersion_})
that the manifold $\tilde V$ is an
algebraic cone over a projective
orbifold $X$. In other words,
there exists an ample bundle
$L$ on a projective orbifold $X$
such $\tilde V$ is biholomorphic to 
the total space $\Tot(L^{\neq 0})$
of the $\C^*$-bundle of all non-zero
vectors in $L$. However, the
$\Z$-action giving $V$ does not
necessarily preserve the algebraic
cone structure. 

The main complex-geometric
invariant of the Vaisman manifold
is the canonical foliation
$\Sigma$ (\ref{_cano_foli_Definition_}),
which is independent from its metric structure
(\cite{_Verbitsky_vanishing_}). It turns out that 
 we can recover the $\Z$-action from 
the vector fields $\langle \theta^\sharp, I(\theta^\sharp\rangle$ 
which trivialize the bundle $T\Sigma$
(\ref{_Lee_field_Remark_}).
The following theorem was proven in
\cite{_OV:flow_}.

\hfill

\theorem\label{_confo_then_Z_Theorem_}
Let $V$ be a compact LCK manifold, and
$\rho:\; \R \times V \arrow V$ a holomorphic conformal flow
of diffeomorphisms on $V$. Assume that $\rho$
is lifted to a flow of non-isometric homotheties
on the smallest K\"ahler covering $\tilde V$ of $V$.
Let $\Gamma\subset \Aut (\tilde V)$ be the deck transform
group (that is, the monodromy), $G$ the
$C^0$-closure of the group generated by $\rho$
in $\Diff(V)$, and $\tilde G$ its preimage in
$\Diff(\tilde V)$. Then $\tilde G$ is connected
and contains $\Gamma$.

\proof 
\cite[Theorem 2.1]{_OV:flow_}. \endproof

\hfill

Let $\theta^\sharp$ be the Lee field, that is, the vector field
dual to the Lee form, $G$ the $C^0$-closure of the flow of $\theta^\sharp$
and $I(\theta^\sharp)$ in the
group of diffeomorphisms of $V$. Since the canonical foliation is unique, $G$ only depends on the complex structure of $V$. Let  $\tilde G$ be the lift of $G$ 
to the $\Z$-cover $\tilde V$. Since $\theta^\sharp$ is Killing,
we can apply \ref{_confo_then_Z_Theorem_}
and obtain that $\Gamma\subset \tilde G$.

\hfill

The elliptic curves on a Vaisman manifold $V$ 
coincide with the closed 2-dimensional orbits of $G$.
We are going to count the number of such orbits.

For any subgroup $\Z=\Gamma' \subset \tilde G$ 
sufficiently close to $\Gamma$, the
quotient $\tilde V/\Gamma'$ is also Vaisman;
however, as shown in \cite[Theorem 4.5]{_Ornea_Verbitsky_immersion_},
there is a dense set of subgroups
$\Z=\Gamma' \subset \tilde G$ 
such that the canonical foliation 
on the quotient $V':=\tilde V/\Gamma'$
is quasi-regular, that is, it has compact leaves. 
In this case the leaf space of the canonical
foliation is a projective orbifold, and the
Lee and anti-Lee fields are contained in 
$\tilde G$.

Clearly, the number of 2-dimensional orbits of $G$ on $V$,
of $\tilde G$ on $\tilde V'$ and of $G':= \tilde G/\Gamma$
on $V'$ is the same. However, each $k$-dimensional
orbit of $G'$ on $V'$ is elliptically fibered
over a $(k-2)$-dimensional orbit of $G'$ on
the leaf space $X:=V'/\Sigma_{V'}$ of the
corresponding canonical foliation. Therefore,
the number of elliptic curves on $V$
is the same as the number of fixed points
of $G'$ acting on $X$. This number can
be found using the
Bia\l ynicki-Birula's and Fontanari's theorems
as in \ref{_Fontanari_Theorem_}, bringing us the
following result.

\hfill

\theorem\label{_Vaisman_number_Theorem_}
Let $V$ be a Vaisman manifold, $V'$
its quasi-regular deformation, and $X:=V'/\Sigma_{V'}$ 
the corresponding K\"ahler orbifold. Assume that
$V$ admits only finitely many elliptic curves.
Then the number of elliptic curves on $V$ 
is equal to $\sum_{i=0}^{\dim_\R X} b_i(X)$.
\endproof

\hfill

Using the Lefschetz $\goth{sl}(2)$-action on the cohomology
of $X$, we obtain a convenient lower bound.

\hfill

\corollary\label{_Vaisman_bound_from_Lfsch_Corollary_}
Let $V$ be a Vaisman manifold which contains precisely 
$r$ elliptic curves. Then $r\geq \dim_\C V$.
\endproof

\hfill

\noindent{\bf Acknowledgment:} We are grateful to Gianluca Bande and Charles Boyer  
for drawing our attention to the reference \cite{_Rukimbira_} and to Oliver Goertsches and Hiraku Nozawa for clarifying us a part of their work \cite{_Goertsches_Nozawa_Tobin_}.

\hfill

{\scriptsize

\hfill

{\small
	
	\noindent {\sc Liviu Ornea\\
		University of Bucharest, Faculty of Mathematics, \\14
		Academiei str., 70109 Bucharest, Romania}, and:\\
	{\sc Institute of Mathematics "Simion Stoilow" of the Romanian
		Academy,\\
		21, Calea Grivitei Str.
		010702-Bucharest, Romania\\
		\tt lornea@fmi.unibuc.ro, \ \  liviu.ornea@imar.ro}
	
	\hfill

	\noindent {\sc Misha Verbitsky\\
		{\sc Instituto Nacional de Matem\'atica Pura e
			Aplicada (IMPA) \\ Estrada Dona Castorina, 110\\
			Jardim Bot\^anico, CEP 22460-320\\
			Rio de Janeiro, RJ - Brasil }\\
		also:\\
		Laboratory of Algebraic Geometry, \\
		Faculty of Mathematics, National Research University 
		Higher School of Economics,
		6 Usacheva Str. Moscow, Russia}\\
	\tt verbit@verbit.ru, verbit@impa.br }
}

\end{document}